\documentclass[graybox]{svmult}

\usepackage{mathptmx}       
\usepackage{helvet}         
\usepackage{courier}        
\usepackage{type1cm}        
\usepackage{makeidx}         
\usepackage{graphicx}        
\usepackage{multicol}        
\usepackage[bottom]{footmisc}
\usepackage{amssymb, amsmath, amsfonts, latexsym, enumerate}
\usepackage[utf8]{inputenc} 
\def \sur#1#2{\mathrel{\mathop{\kern 0pt#1}\limits^{#2}}}
\def \el{\sur{=}{(law)}}
\makeindex 

\begin{document}
\title*{Large deviations for clocks of  self-similar 
 processes}
\author{Nizar Demni, Alain Rouault and Marguerite Zani}
\institute{Nizar Demni \at Institut de Recherche Mathematiques de Rennes, Université Rennes 1, Campus de Beaulieu 35042 Rennes, France; \email{Nizar.Demni@univ-rennes1.fr}
\and Alain Rouault \at Université Versailles-Saint-Quentin, LMV  UMR 8100, Bâtiment Fermat,  45 Avenue des Etats-Unis, 78035 Versailles-Cedex, France; \email{Alain.Rouault@uvsq.fr}
\and Marguerite Zani \at Université d'Orléans, UFR Sciences,
Bâtiment de mathématiques - Rue de Chartres
B.P. 6759 - 45067 Orléans cedex 2, 
France; \email{marguerite.zani@univ-orleans.fr}}
\maketitle
\centerline{\textit{In memoriam, Marc Yor}}
\bigskip
\abstract{
The Lamperti correspondence gives a prominent role to two random time changes: the exponential functional of a Lévy process drifting to $\infty$ and its inverse, the  clock of the corresponding positive 
self-similar
 process. We describe here asymptotical properties of these clocks 
in large time, extending the results of \cite{zani}.}
\section{Introduction}
This problem is an extension of a question  raised by Marc Yor during the defense of the thesis of Marguerite, under the supervision of Alain, long time ago in 2000. 
The last part of this thesis was dedicated to the study of large deviations principles (LDP) for Maximum Likelihood Estimates of diffusion coefficients (for squared-radial Ornstein--Uhlenbeck processes, squared Bessel processes and Jacobi processes).  
The main tool there was a convenient Girsanov change of probability.
This method allowed to convert 
the computation of Laplace transform of some  additive functionals into the computation of Laplace transform of  a single variable. This trick was used before  in  \cite{yorzurich1} page 26, or \cite{manyor}, page 30, where Marc Yor called it "reduction method". 
In those times, Marc was interested in exponential functionals of Brownian motion, Lamperti transform and Asian options and he
guessed that this LDP could be applied to the Bessel clock and solved effectively the problem with Marguerite in \cite{zani} a couple of weeks later. Marc suggested to extend it to the Cauchy clock and gave a sketch of proof, but a technical difficulty stopped  the project. 
Meanwhile, Marguerite and Nizar published a paper \cite{demni_zani} on Jacobi diffusions where the reduction method  is again crucial.
Recently, the three of us  felt the need to revisit the problem of Cauchy clock with the hope of new ideas. 

We planned to discuss with Marc, and promised him we would keep him informed of our progress. We did not have the time... 
\par

In this paper, we extend the methods and results of \cite{zani} to a large class of clocks issued from positive self-similar Markov processes.
In  Section 2 we recall 
 some basic results about the Lamperti correspondence between these processes and Lévy processes
 and we give the definition of the clocks. 
We also define a generalized Ornstein-Uhlenbeck process which will be useful in the sequel. In Section 3 and Section 4, we  present the main results: Law of Large Numbers and Large Deviations for the clocks.
 In Section 5 we show some examples illustrating our main theorems.

\section{Positive self-similar Markov 
 processes and Lamperti transformation}

In \cite{lamp} Lamperti defined (positive) semi-stable process which are nowadays called  positive self-similar 
 Markov process.

\begin{definition}
For $\alpha > 0$, a 
 positive self-similar 
 Markov process (pssMp) 
  of index $\alpha$,   is a $[0, \infty)$-valued strong Markov process $(X, \mathbb Q_a), a > 0$ with càdlàg paths, fulfilling the scaling property
\begin{equation}
\label{self}
\left(\{b X_{b^{-\alpha}t}, t \geq 0\}, \mathbb Q_a\right) \el \left(\{X_t , t \geq 0\}, \mathbb Q_{ba}\right) 
\end{equation}
for every $a, b > 0$. 
\end{definition}

Lamperti \cite{lamp} has shown that these processes can be connected  to Lévy processes by a one-to-one correspondence, that we develop below.
We refer to 
 Kyprianou \cite{Kyp} especially Chapter 13 for properties of Lévy processes and pssMp. One can also see \cite{bertoin2004exponential} for the Lamperti's correspondence. One can notice that there is a little confusion in the notion of index of these processes. In \cite{lamp}, \cite{bertoin2002entrance} and \cite{bertcaba}, the index is $1/\alpha$, and in \cite {bertoin2004exponential}, \cite{cab} and \cite{Kyp}, the index is $\alpha$. We take this latter convention.
 These processes have a natural application in the theory of self-similar fragmentations : see \cite{bertfrag}, references therein, and \cite{cr}. For other areas of application, such as diffusions in random environments, see Section 6 of \cite{bertoin2004exponential}.
\subsection{From $X$ to $\xi$}
Any pssMp $X$   which never reaches the boundary state $0$
 may be expressed as the
exponential of a Lévy process  not drifting to $-\infty$, time changed by the inverse of its exponential functional. 
More formally, if $(X, (\mathbb Q_a)_{a > 0})$ is  a  pssMp of index $\alpha$  which never reaches $0$,  set
\begin{equation}
T^{(X)}(t) = \int_0^t \frac{ds}{X_s^\alpha} \  , \ (t \geq 0)
\end{equation}
and let $A^{(X)}$ be  its inverse, defined by
\begin{equation}
A^{(X)}(t) = \inf \{u \geq 0 : T^{(X)}(u) \geq t \}\,,
\end{equation}
and let $\xi$ be the process defined by
\begin{equation}
\xi_t = \log X_{A^{(X)}(t)} - \log X_0\ , \ (t\geq 0)\,.
\end{equation}
Then, for every $a >0$, the distribution of $(\xi_t, t \geq 0)$ under $\mathbb Q_a$ does not depend on $a$ and is the distribution of a Lévy process starting from $0$.

Moreover, if we set
\begin{equation}
\label{expfunc}
\mathcal{A}^{(\xi)} (t) = \int_0^t e^{\alpha \xi_s}ds
\end{equation}
and $\tau^{(\xi)}$ its inverse defined by
\begin{equation}
\label{invexpfunc}
\tau^{(\xi)}(t) = \inf\{u \geq 0 : \mathcal{A}^{(\xi)}(u) \geq t\}
\end{equation}
we have
\begin{equation}
\label{fund}
\tau^{(\xi)} (t) = T^{(X)} (t X_0^\alpha)\,.
\end{equation}
Let us remark that the self-similarity property (\ref{self}) leads to the following relation: 
\begin{equation}
\label{invtau}
\hbox{the law of} \ \ T^{(X)} (.\  b^{-\alpha}) \ \ \hbox{under} \ \  \mathbb Q_a\ \ \hbox{is}\ \hbox{the law of} \ \ T^{(X)}(.) \ \ \hbox{under} \ \ \mathbb Q_{ba}\,.
\end{equation}
\subsection{From $\xi$ to $X$}
\label{xiX}
Let  $\xi_t$ be a  Lévy process starting from $0$ and let $\mathbb P$ and $\mathbb E$ denote the underlying probability and expectation, respectively. 
 Fix $\alpha > 0$. Let $\mathcal{A}^{(\xi)}$ be its exponential functional defined by (\ref{expfunc}).
When $\xi$ drifts to $ -\infty$, this functional is very popular in mathematical finance (see \cite{yor2001exponential}), with important properties of the perpetuity $\mathcal{A}^{(\xi)}(\infty)$. Here, we rather assume that  $\xi$ does not drift to $-\infty$ i.e. satisfies,   
$\limsup_{t\uparrow \infty} \xi_t = \infty$. We define the
 inverse process $\tau^{(\xi)}$ of $\mathcal A^{(\xi)}$ by (\ref{invexpfunc}).

For every $a >0$, let $\mathbb Q_a$ be the law under $\mathbb P$ of the time-changed process
\begin{equation}
\label{xitoX1}
X_t = a\exp \xi_{\tau^{(\xi)}(ta^{-\alpha})} \ , \ (t\geq 0)\,,
\end{equation}
then $(X , (\mathbb Q_a)_{a > 0})$ is a  pssMp of index $\alpha$  which never reaches $0$ and we have the fundamental relation
 (\ref{fund}). 
\subsection{Index and starting point}
\label{index}
 If $(X, (\mathbb Q_a)_{a >0})$ is  an pssMa of index $\alpha$, then  the process $Y =( X^{\alpha}, (\mathbb Q_{a^\alpha})_{a >0})$, is a pssMp of index $1$.
 Conversely if $(Y, (\mathbb Q_a)_{a > 0})$ is a pssMp if index $1$
then, for any $\alpha > 0$, 
the process 
$(X = Y^{1 / \alpha} ,(\mathbb Q_{a^{1/\alpha}})_{a > 0})$
 is a pssMp of index $\alpha$. Hence \underline{we may and will assume 
without loss of generality, that $\alpha = 1$}, except some examples in Section \ref{five}.

Let us stress also the fact that the distribution of  $T^{(X)}$ 
is considered under $\mathbb Q_a$ for $a >0$, and the distribution of $\tau^{(\xi)}$ is considered under $\mathbb P$.

\subsection{Ornstein-Uhlenbeck process}
Having defined  the above pair of associated processes $(\xi, X)$ it is usual to consider a third process, very useful in the sequel. To introduce it, we need additional properties on the pssMp $X$. 
In Theorem 1 of \cite{bertoin2002entrance}, Bertoin and Yor have shown that whenever the support of $\xi$ is not arithmetic and $\mathbb E(\xi_1)>0$, then as $a\downarrow 0$, the probability measures $\mathbb Q_a$ (defined in section \ref{xiX}) converge in the sense of finitely dimensional marginals to a probability measure denoted by $\mathbb Q_{0}$. This latter measure is an entrance law for the
 semigroup $Q_t f(x) = \mathbb E_x f(X_t)$ %
and satisfies
\begin{equation}
\label{magicI}
\mathbb Q_0(f(X_t)) = \frac{1}{\mathbb E \xi_1} \mathbb E[I_\infty^{-1}f(t/I_\infty)]
\end{equation}
where 
\[I_\infty 
 = \int_0^\infty \E^{-\xi_s} \D s\,.\]
Now we define the associated generalized Ornstein-Uhlenbeck (OU) process by
\[U(t) := \E^{- t}X(\E^{t})\,.\] 
This process has been studied in \cite{CPY2}, and further on by \cite{Riv} and \cite{cr}.
It is  strictly stationary, Markovian and ergodic under $\mathbb E_0$,  and its invariant measure is the law of $X_1$ under $\mathbb Q_0$, i.e.
\[\mu(f) =  \frac{1}{\mathbb E \xi_1} \mathbb E[I_\infty^{-1}f(I_\infty)]\,.\]
\medskip

The purpose of this paper is the study of the asymptotic behavior of the  process $T^{(X)}$, called the clock of the pssMp $X$, as $t \rightarrow \infty$  under $\mathbb Q_a$ for any $a >0$, or equivalently  the asymptotic behavior of the  process 
$\tau^{(\xi)}$ , as $t \rightarrow \infty$  under $\mathbb P$.  
 The two points of view are complementary:  the clock as a functional of the pssMp or the clock as a functional of the Lévy process.
We will restrict us to the case of $\mathbb E \xi_1 \in (0, \infty)$, which ensures a Law of Large Numbers (LLN) and a Large Deviation Principle (LDP). 

\section{Law of Large Numbers}
The first step of our study is the a.s. convergence of our functionals.
\begin{theorem}
\label{thmLLN}
If $\mathbb E \xi_1$ is finite positive (hence $= \psi'(0) > 0)$,
then for every $a >0$, we have 
\begin{equation}
\label{LLNX}
\frac{T^{(X)}(t)}{\log t} \rightarrow (\mathbb E \xi_1)^{-1} \,, \ \ \mathbb Q_a-\hbox{a.s.}
\end{equation}
or equivalently
\begin{equation}
\label{LLN}
\frac{\tau^{(\xi)}(t)}{\log t} \rightarrow (\mathbb E \xi_1)^{-1} \,, \ \ \mathbb P-\hbox{a.s.}
\end{equation}%
\end{theorem}

\proof From (\ref{fund}) it is enough to prove (\ref{LLNX}). 
Let us now repeat verbatim the trick from 
\cite{cr} Section 2.
By the ergodic theorem, we have that for $f \in L^1(\mu)$:
\[\frac{1}{t}\int_0^t f(U_s) \D s \rightarrow \mu(f)  \,, \ \ \mathbb Q_0-\hbox{a.s.}\]
or, by scaling
\begin{equation} 
\frac{1}{\log t}\int_1^t f(u^{-1}X_u) \frac{\D u}{u} \rightarrow \mu(f) \,, \ \ \mathbb Q_0-\hbox{a.s.}
\end{equation}
For $f(x) = x^{-1}$ this yields
\[\frac{1}{\log(1+t)} \int_1^{1+t} \frac{\D u}{X_u} \rightarrow (\mathbb E\xi_1)^{-1} \ \,,\ \ \mathbb Q_0-\hbox{a.s.}\]
If $\mathcal H$ denotes the set where the convergence holds, the Markov property yields:
\[\mathbb E_0\left[\mathbb Q_{X_1}\left(\frac{T^{(X)}(t)}{\log(1+t)} \rightarrow \  (\mathbb E\xi_1)^{-1} \right)\right] 
= \mathbb P_0(\mathcal H)=1\,.\]
If we remember that $\mu$ is the law of $X_1$ under $\mathbb Q_0$, we have 
\[\mathbb Q_x \left(\frac{T^{(X)}(t)}{\log(1+t)}
 \rightarrow \  (\mathbb E\xi_1)^{-1} \right) = 1\,, \ \ \mu-\hbox{a.e.}\,.\]
 Fixing such an $x$ and using the scaling property (\ref{invtau}), we get  (\ref{LLNX}).
\qed
\medskip

It is then natural  to look for a Large Deviation Principle (LDP) to characterize the speed of this convergence.

\section{LDP}
We present a LDP for the clock $T^{(X)}(t)$  (Theorem \ref{thmLDP}), whose
proof relies on the reduction method (Lemma \ref{ncgf}). The remaining of the section is devoted to some
illustrative remarks and the proof of a technical lemma. In view of (\ref{fund}), the following statements are valid when $T^{(X)}$ and $\mathbb Q_a$ are replaced by  $\tau^{(\xi)}$ and $\mathbb P$, respectively.

\subsection{Main result}

\begin{theorem}
\label{thmLDP}
Assume\footnote{Domains are taken open par convention} that dom $\psi = (m_-, m_+)\ni 0$ and that $\psi'(0) > 0$. Set
\[m_0 = \inf\{ \theta : \psi' (\theta) > 0\}\ \ , 
\tau_+ =  \frac{1}{\psi'(m_+)} \ , \ \tau_0 =  \frac{1}{\psi'(m_0)} \ \hbox{and} \ \Delta = (\tau_+, \tau_0)\,,\]
where $1/ \psi'(\pm \infty) := \lim_{m\rightarrow \pm \infty} m/\psi(m)$.
Let $a >0$.
\begin{enumerate}
\item
Then, for every $x \in  \bar\Delta$
\begin{equation}
\label{weakLDP}
\lim_{\varepsilon \rightarrow 0}\lim_{t\rightarrow \infty} \frac{1}{\log t} \log \mathbb Q_a \left(\frac{T^{(X)}(t)}{\log t} \in [x-\varepsilon, x + \varepsilon]\right) = - \mathcal I(x)
\end{equation} where 
\begin{equation}\label{defI}\mathcal I(x) = \sup_{m\in (m_0, m_+)} \{m- x\psi(m)\}\,.\end{equation}
\item Moreover 
\begin{equation}
\lim_{A\rightarrow \infty} \limsup_{t\rightarrow \infty}\frac{1}{\log t} \log \mathbb Q_a \left(\frac{T^{(X)}(t)}{\log t} > A\right) = - \infty\,. 
\end{equation}

\item
If either $\Delta = (0, \infty)$ or the complement of $\bar \Delta$ is exponentially negligible, i.e. if
\begin{equation}
\label{assum}
\forall \varepsilon > 0 \ \ \  \limsup_{t\rightarrow \infty}\frac{1}{\log t} \log \mathbb Q_a \left(\inf_{ x \in \bar\Delta} \left\{ \left|x- \frac{T^{(X)}(t)}{\log t}\right| \right\} > \varepsilon\right) = - \infty
\end{equation} 
then 
 the family of distributions of $T^{(X)}(t)/ \log t$ under $\mathbb Q_a$ satisfies the LDP on $[0, \infty)$ at scale $\log t$ with good rate function :
\[\widetilde{\mathcal I}(x) =
\begin{cases}
\mathcal I(x)\;\;\;\;&\mbox{if}\; x \in \bar\Delta,\\
\infty\;\;&\mbox{ otherwise.} 
\end{cases}\]
\end{enumerate}
\end{theorem}

The following proposition describes some properties 
 of  $\mathcal I$ which are direct consequences of \eqref{defI}. The proof is left to the reader.

\begin{proposition}
\label{classification}
\begin{enumerate}
\item
 If $\psi^*$ denotes the Fenchel-Legendre dual of $\psi$ defined by
\[\psi^*(x) = \sup_{m\in \mathbb R} {mx -\psi(m)} \ \ , \   (x\in \mathbb R)\,,\] then
\begin{equation}
\label{pair}
{\mathcal I}(x) = x\psi^*(x^{-1}) \ \ ,  \ (x \in \Delta)\,.\end{equation}
\item Let $\tau_e = 1/ \mathbb E \xi_1 \in (\tau_+, \tau_0)$. The function $\mathcal I$ is convex on $\Delta$, decreasing on 
$(\tau_+, \tau_e)$
, has a minimum $0$ at $\tau_e$
, and is increasing on 
$(\tau_e, \tau_0)$.
\item
\begin{enumerate}
\item
If $m_0$ is a true minimum, i.e. if $m_0 > -\infty$ and $\psi'(m_0) = 0$, then $\tau_0= \infty$, and as $x \rightarrow \infty$, $\mathcal I$ admits 
   the asymptote $y= -x \psi(m_0) + m_0$. 
\item
If $m_0 = -\infty$
and $\psi'(m_0) \in (0, \infty)$ 
, set $\lim_{m\rightarrow m_0} \psi(m)- m\psi'(m_0) = -b\in [-\infty, 0)$. Then $\tau_0 < \infty$ and 
$\mathcal I(\tau_0) = b\tau_0$ and $\mathcal I'(\tau_0) = - \psi(m_0)= \infty$.
\item
If $m_0 = -\infty$ and $-\infty < \psi(m_0) < 0$, hence $\psi'(m_0) = 0$, then $\tau_0 = \infty$ and $\mathcal I(x) / x \rightarrow -\psi(m_0)$ and $\mathcal I(x) +\psi(m_0) x \rightarrow -\infty$, as $x \rightarrow \infty$.
\end{enumerate}
\item
\begin{enumerate}
\item
If $m_+ < \infty$ and $\psi(m_+) = \infty$, then $\psi'(m_+) = \infty$. We get $\tau_+ =0$ and $\mathcal I(0) = m_+$ with $\mathcal I'(0) = \infty$.
\item
If $m_+ = \infty$ and $\psi'(m_+) < \infty$, set $\lim_{m \rightarrow m_+} \psi(m) - m \psi'(m_+) = - b\in [-\infty, 0)$. Then 
$\tau_+ > 0$, $\mathcal I(\tau_+) = b\tau_+ \leq \infty$ and $\mathcal I'(\tau_+) = -\psi(m_+) = - \infty$. 
\item If $m_+ = \infty$, with $\psi(m_+) = \infty$ and $\psi'(m_+) = \infty$ then $\tau_+ = 0$ and 
$\mathcal I(\tau_+) = m_+ =  \infty$.
\end{enumerate}
\end{enumerate}
\end{proposition}

\proof {\it of Theorem \ref{thmLDP}}

1) A slight adaptation of the G\"{a}rtner-Ellis method (\cite{demboz98} Th. 2.3.6) allows to deduce (\ref{weakLDP}) from the asymptotic behaviour of the normalized log-Laplace transform of $T^{(X)}(t)$, given by the following lemma. Recall that $\psi$ is increasing on $(m_0, m_+)$.
\begin{lemma}
\label{ncgf} 
For $\theta \in \left(-\psi(m_+), - \psi(m_0)\right)$, set $L_t(\theta) = \log \mathbb E_a  \exp(\theta T^{(X)}(t))$. Then as $t \rightarrow \infty$ we have
\begin{equation}
\frac{1}{\log t} L_t (\theta) \rightarrow L(\theta)
\end{equation}
where 
\begin{equation}
\label{Lpsi}
L(\theta) = -m \Longleftrightarrow \theta = - \psi(m)\,.
\end{equation}
\end{lemma}
The function $L$ is differentiable on $\left(-\psi(m_+), - \psi(m_0)\right)$ and satisfies $L'(\theta) = 1/\psi'(m)$, so that the range of $L'$ is precisely $\Delta$. Then the left-hand-side of (\ref{weakLDP}) admits the limit  
$\mathcal I(x) = x\theta - L(\theta)$ where $\theta$ is the unique solution of $L'(\theta) = x$ i.e., thanks to (\ref{Lpsi}), 
$\mathcal I(x) = -x\psi(m)+ m$ which is exactly the right-hand-side of (\ref{weakLDP}).

\proof {\it of Lemma \ref{ncgf}}

To compute the Laplace transform of $T^{(X)}(t)$, we use a Girsanov type change of probability.

For $m \in (m_-, m_+)$ let 
\[\psi_m (\theta) = \psi(m+\theta) - \psi(m)\,,\]
and let $\{(\xi_t , t \geq 0) ; \mathbb P^{(m)}\}$ be a Lévy process starting from $0$ whose exponent is $\psi_m$ 
(Esscher transform).
Finally let   $\{(X_t , t \geq 0) ; (\mathbb Q^m_a)_{a > 0}\}$  be  the associated pssMp.

Besides, from  \cite{CPY} the following relation between the infinitesimal generators $L^{\xi}$ of $\xi$ and $L^X$ of $X$ under $\mathbb Q_a$:

\[L^Xf(x)=\frac{1}{x}L^{\xi}(f\circ \exp)(\log x)\,,\]
implies (see   (2.7) in \cite{CPY}),

\[L^X f_m = \psi(m) f_{m-1}\]
for $f_m(x) = x^m$ so that
\begin{eqnarray}
\mathbb Q^m_a|_{\mathcal F_t}= \left(\frac{X_t}{a}\right)^m \exp \left(- \psi(m)\int_0^t \frac{\D s}{X_s}\right)\  .\ 
\mathbb Q_a |_{\mathcal F_t}\ \ , \  (t \geq 0)\,. 
\end{eqnarray}
 We deduce that 

\begin{equation}
\mathbb E_a \exp  \left(- \psi(m)\int_0^t \frac{ds}{X_s}\right)  = 
 a^m \mathbb E_a^m (X_t)^{-m}
\end{equation}
and, owing to the scaling property
\begin{equation}
\label{magic}
\mathbb E_a \exp  \left(- \psi(m)\int_0^t \frac{ds}{X_s}\right)  = 
 a^m t^{-m}\mathbb E^m_{a/t} (X_1)^{-m}\,.
\end{equation}
Let us choose $m$ such that $\psi'(m) > 0$.

We can know use the results on entrance boundary, for example  \cite{bertoin2002entrance}. Since the new Lévy process  $\xi^{(m)}$ has a positive mean
\[\mathbb E \xi^{(m)} = \psi'_m(0) = \psi'(m) >0\]
then Theorem 1 i) therein entails that
\begin{equation}
\label{limentrance}\lim_{t \uparrow \infty} \mathbb E^m_{a/t} (X_1)^{-m}= \mathbb E^m_0  (X_1)^{-m}\,. \end{equation}

\begin{lemma}
\label{finite}
With the notations of Theorem \ref{thmLDP},  the quantity $F(m) :=  \mathbb E^m_0  (X_1)^{-m}$ is finite for $m \in (m_0, m^+)$.
\end{lemma}

From (\ref{magic})
\[L_t(\theta) = -m \log t  + m \log a +\log \mathbb E^m_{a/t} (X_1)^{-m}\]
 and then, as $t \rightarrow \infty$, thanks to Lemma \ref{finite}:
\begin{equation}
\label{defL}\frac{L_t (\theta)}{\log t} \rightarrow L(\theta) = -m\,.\end{equation}
This ends the proof of Lemma \ref{ncgf} and consequently the proof of part 1) in Theorem \ref{thmLDP}, up to the result of Lemma \ref{finite} whose proof is postponed.
\medskip

2) The statement   is a consequence of the Chernov inequality. Indeed 
for fixed    $\theta \in (0, -\psi(m_0))$
\[\log \mathbb Q_a\left(\frac{T^{(X)}(t)}{\log t} > - A\right) \leq -\theta A + L(\theta)\]
for every  $\theta \in (0, -\psi(m_0))$.
\medskip

3) Let us assume that  $\Delta = (0, \infty)$. To get a full LDP, it is enough to prove that
\[\lim_{A \rightarrow \infty} \lim_{t\rightarrow \infty} \frac{1}{\log t} \log \mathbb P\left(\frac{\tau(t)}{\log t} > A\right) = - \infty\,,\]
(exponential tightness) but it is a consequence of the Chernov inequality:
\[\log \mathbb Q_a\left(\frac{T^{(X)}(t)}{\log t} > - A\right) \leq -\theta A + L(\theta)\]
for every  $\theta \in (0, -\psi(m_0))$. 

Assume now that the condition (\ref{assum}) is satisfied. For every $x$ in the complement of $\bar\Delta$ we have
\[\lim_{\delta \rightarrow 0} \lim_{t\rightarrow \infty} \frac{1}{\log t} \log \mathbb Q_a\left(\frac{T^{(X)}(t)}{\log t} \in [x-\delta, x+\delta]\right) = -\infty\,, \]
so that a weak LDP is satisfied on $[0, \infty)$. To strenghten it is enough to apply the above reasoning.
\qed

\subsection{Remarks}
\begin{remark}(Reciprocal pairs)
The   relation (\ref{pair}) between two rate functions is known to hold for pairs of inverse processes. This problem arose historically when people deduced a LDP for renewal processes whose interrival times have sums satisfying a LDP. It was extended to more general processes, see \cite{DW} and the bibliography therein. 

If $(\xi_t)$ is a spectrally negative Lévy process   with Laplace exponent $\psi$, then $t^{-1}\xi_t$  satisfies the LDP at scale $t$ with rate function $\psi^*$. The subordinator 
\begin{equation}\label{taut}\hat\tau(u) = \inf\{ t > 0 : \xi_t > u\}\end{equation} has a Laplace exponent which is exactly $L$ defined by he relation (\ref{Lpsi})
\begin{equation}
\label{recip}\psi(m) = -\theta \ , \ L(\theta) = -m\,.\end{equation}
Then $t^{-1}\hat\tau(t)$ satisfies the LDP at scala $t$ with rate function  $x\to x\psi^*\left(\frac{1}{x}\right)$. This can be seen    as an application of Theorem 1 of \cite{DW}).

The relation (\ref{invexpfunc}) between $\mathcal A^{(\xi)}_u$ and $\tau^{(\xi)}(t)$:
\[\tau^{(\xi)}(t)=\inf\{u\geq 0: \mathcal A^{(\xi)}_u=\int_0^ue^{\xi_s}ds\geq t\}\ \ ,  \ (t\geq 0)\]
\noindent
 is more involved than \eqref{taut}. An alternative proof of our Theorem \ref{thmLDP} about LDP for $(\log t)^{-1}\tau^{(\xi)}(t)$, in the spirit of   Theorem 1 of \cite{DW}), would require an LDP for  $t^{-1}\log \mathcal A_t^{(\xi)}$ which is far from obvious. At first glance, the lower bound of large deviations seems accessible since we have for convenient $x$ and $\delta$,
\[\left(\xi_s/s \in [x-\delta, x+\delta] \ \forall s \in [0,u]\right) \Rightarrow \left(\mathcal A^{(\xi)}_u \in \left[\frac{\E^{(x-\delta)u} -1}{(x-\delta)} ,\frac{\E^{(x+\delta)u} -1}{(x+\delta)}\right]\right)\,.\]

Conversely, we can apply Theorem 1 of \cite{DW} to get 
\begin{proposition}
The family of distributions of $t^{-1}\log \mathcal A^{(\xi)}_t$ (under $\mathbb P$) satisfies the LDP at scale $t$, with good rate function $\psi^*$.
\end{proposition}

\noindent To end this first remark, let us add two more comments:
\begin{itemize}
\item[$\quad\quad\bullet$] \hskip15ptThe statement of Lemma \ref{ncgf}, viewed in terms of $\xi$ can be rephrased as   the striking relation
\[\lim_{t\rightarrow \infty} \frac{1}{\log t} \mathbb E\exp\theta \tau^{(\xi)}(t) = \log \mathbb E \exp \theta \hat\tau(1)\,.\]
\item[$\quad\quad\bullet$] \hskip15pt
The relation (\ref{recip}) is also well known in the study of exponential families (inverse families or reciprocal pairs) , see  \cite{kuchler1997exponential} Section 5.4 and  \cite{LetacMora} Section 5A.
\end{itemize}
\end{remark}
\begin{remark}(Another approach for the LDP)
We saw in a previous section that the core of the proof of the LLN is the convergence of the occupation measure
\[\frac{1}{t} \int_0^t \delta_{U_s} \D s \rightarrow \mu\]
followed by a scaling and by the choose of $f(x) = x^{-1}$ as a test function. It could be natural to look for a LDP according the same line of reasoning. Such a proof would consist of three steps.  First establish an LDP for the law of the occupation measure under $\mathbb Q_0$. The rate functional is expressed  with the help of the infinitesimal generator of $U$. Then apply the contraction by the (non-continous) mapping $\nu \mapsto \nu(f)$, hence solve a variational problem. This would give an LDP for the law of
\[\frac{1}{t} \int_0^t \frac{\D s}{U_s} = \frac{1}{t}\int_1^{e^t} \frac{\D s}{X_s}\,,\]
under $\mathbb Q_0$.
It would remain to convert it into an LDP under $\mathbb Q_a$. In the Appendix of \cite{zani}, the authors give the first rate functional and solved the variational problem, but with the lack of justification for the contraction principle. The result fits with the rate function of Theorem 2.
\end{remark}

\begin{remark}
\label{FLDP}
(Functional LDP)
In the Bessel clock, a functional LDP was stated (\cite{zani} Th. 4.1). Here it is possible to consider the same problem, i.e. the study of the LDP for the sequence of processes
 \[u\in [0,1] \mapsto \frac{1}{n} T^{(X)}(\E^{nu})\]
under $\mathbb Q_a$. The rate function involves an action functional built on $\mathcal I$ and a
 recession function which will be $x \mapsto rx$ where $r =   -\psi(m_0)$ (when it is finite).
\end{remark}

\begin{remark}(Central Limit Theorem)
It is known for the  Bessel clock (Theorem 1.1 in \cite{zani}) but seems unknown otherwise. We can conjecture that
\begin{equation}
\label{conjCLT}
\sqrt{\log t}\left(\frac{T^{(X)}(t)}{\log t} - \frac{1}{\mathbb E \xi_1}\right)\Longrightarrow \mathcal N(0 ; \psi''(0)/ \psi'(0)^3)
\end{equation}
as soon as $\mathbb E (\xi_1)^2 < \infty$ (or $\psi''(0) < \infty$). 
\end{remark}
\subsection{Proof of Lemma \ref{finite}}
From (\ref{magicI}) applied to the Esscher transform we have
 for  $m > m_0$
\[F(m) = \mathbb E^{(m)} (I_\infty^{m-1})\,.\]
The finiteness of $F(m)$ is a consequence of the following Lemma summarizing properties of the moments of the exponential functionals of Lévy processes which we detail for completeness.

\begin{lemma}
\label{zwart}
 Let $\zeta$ be a Lévy process of Laplace exponent $\varphi$ given by
\[\mathbb E \exp(\lambda \zeta_t) = \exp(t\varphi(\lambda))\]
such that $\varphi'(0) \in (0, \infty)$ and let
\[I_{\infty}=\int_0^{\infty}e^{-\zeta_s}ds\,.\] Then
\begin{itemize}
\item[1)]
\begin{equation}
\label{maul}
\mathbb E I_{\infty}^s <\infty \ \hbox{for all} \ s \in [-1,0] \ \ \hbox{and all} \ s > 0 : \varphi(-s) < 0\,.
\end{equation}
\item[2)] For $r>0$, if $\varphi(r)<\infty$ and $\mathbb E I_{\infty}^{-r}<\infty$, then 
\begin{equation}\label{recur}
\mathbb E  I_\infty^{-r-1} = \frac{\varphi(r)}{r}\ \mathbb E  I_\infty^{-r}\,.
\end{equation}
\end{itemize}
\end{lemma}
\proof {\it of Lemma \ref{zwart}}

The statement 1) comes from
 a rephrased part of  a lemma in Maulik and Zwart on the existence of moments of exponential functionals
(see \cite{maulik} Lemma 2.1).
\par\noindent
The statement 2) comes from a  a recursive argument 
due to  Bertoin and Yor (see \cite{bertoin2004exponential}, Theorem 2 (i) and Theorem 3), which we detail here for the sake of completeness. 
Set, for $t \geq 0$
\[J_t = \int_t^\infty \E^{-\zeta_s} \D s \ \ (t \geq 0) \ \  ; \ J_0= I_\infty= \int_0^\infty \E^{-\zeta_s} \D s\]
For all $r> 0$ we have
\begin{equation}\label{IPP1}J_t^{-r}- J_0^{-r} = r\int_0^t \E^{-\zeta_s} J_s^{-r-1} \D s\end{equation}
Besides, from the properties of Lévy processes, we deduce
\[J_s = \E^{-\zeta_s} \hat I(s)\]
and
\[\hat I(s) = \int_0^\infty \E^{-(\zeta_{s+u}-\zeta_s)} \D u \el I_\infty\]
with $\hat I(s)$ independent of $J_s$. 
Plugging into (\ref{IPP1}) we get
\begin{equation}
\label{inj1}
\E^{r\zeta_t} \hat I(t)^{-r} - I_\infty^{-r} = r \int_0^t \E^{r\zeta_s} \hat I(s)^{-r-1} \D s
\end{equation}
Assume that  $\varphi(r) < \infty$ and  $\mathbb E I^{-r-1} = \infty$, taking expectations on both sides of the inequality
\[\E^{r\zeta_t} \hat I(t)^{-r} >  r \int_0^t \E^{r\zeta_s} \hat I(s)^{-r-1} \D s\]
we would get $\mathbb E I^{-r}= \infty$. Consequently, if $\mathbb E I^{-r} < \infty$, we have $\mathbb E I^{-r-1} < \infty$. Moreover, taking   expectations on both sides of  (\ref{inj1}) we get (\ref{recur}).
$\hfill\square$
\medskip\par\noindent
Now, to end the proof of Lemma {\ref{finite}, we apply Lemma \ref{zwart} with $\varphi=\psi_m$.
\par\noindent
For $m>0$, we choose $s=m-1$ hence $\varphi(-s)=\psi(1)-\psi(m)$ which is negative if $m>1$. So from 1) of Lemma \ref{zwart}, $F(m)$ is finite for $m>0$.
\smallskip\par\noindent
Now, if $m\in[m_0,0]$, set $k = \lfloor -m\rfloor$. Then $m+k \in [-1,0]$, and from 1) of Lemma \ref{zwart}, $\mathbb E^{(m)} I_\infty^{m+k} $ is finite.
Besides, we have $\psi(m) < \psi(-j)$ for all  integers $j$ such that $0 \leq j \leq  k$, hence $\varphi(-m-j) = \psi(-j) - \psi(m) > 0$ and applying (\ref{recur}) recursively, we have
\begin{equation}
F(m) = \frac{\varphi(-m)\varphi(-m-1)  \cdots \varphi(-m-k)}{(-m)(-m-1)\cdots (-m-k)}\mathbb E^{(m)} I_\infty^{m+k} 
\end{equation}
and this quantity is finite.$\hfill\square$

\section{Examples}
\label{five}
For all the following examples, we give the Laplace exponent $\psi$ and the parameters $\tau_+$, $\tau_e$. In  Sections  \ref{Bessel}, and \ref{Poisson} we obtain an explicit expression for $\tau_0$ and $\mathcal I$. In Section \ref{threemore} we do not have an explicit expression for $\mathcal I$ (and sometimes for $\tau_0$). For all examples (but one) we identify the behaviour at the boundary, according to the classification of Proposition \ref{classification}.

\subsection{The Brownian process with drift  (Bessel clock \cite{zani})}
\label{ex1}
\label{Bessel}

For $\nu > 0$, we consider the Lévy process \[\xi_t = 2 B_t + 2\nu t\,.\]
The pssMp 
 $(X_t)$  is the squared Bessel process of dimension $d= 2(1 +\nu)$. Its  index is $\alpha = 1$ and it is the only continuous 
pssMp
(of index $1$). The clock is related to the scale function of a diffusion in random environment (see \cite{yor2001exponential} chap. 10). 

We have
\[\psi(m) = 2m(m+\nu) , \ m_\pm = \pm \infty,  \ \psi(m_\pm) = +\infty\,.\]
Here $m_0 = -\nu/2$ and $\psi(m_0) = -\nu^2/4$, so that
\[\tau_+ = 0 , \ \ \tau_0 = \infty, \ \ \tau_e = \frac{1}{2\nu} , \ \ \Delta =(0, \infty)\,.\]
The rate function is
\[\mathcal I(x) = \frac{(1-2\nu x)^2}{8x}\  \ , \ (x > 0)\,.\]
Notice that here the function $L$ may be obtained by explicit inversion
\[L(\theta) = \frac{-\nu + \sqrt{\nu^2 -2\theta}}{2} \ \ , \   (\theta <  \nu^2/2)\,.\]
The  boundary $\tau_+=0$ is in the situation 4c) of Proposition \ref{classification} and the boundary $\tau_0$ is in situation 3a, with asymptote  $y = \frac{\nu^2}{2}x - \frac{\nu}{2}$.
The minimum of $\mathcal I$, reached in $\tau_e$, corresponds to the LLN (Theorem \ref{thmLLN}) which is in 
\cite{RY} Exercise (4.23) Chap. IV and in \cite{zani} Theorem 1.1.
\begin{figure}[ht!]
\begin{center}
\includegraphics[height= 150pt]{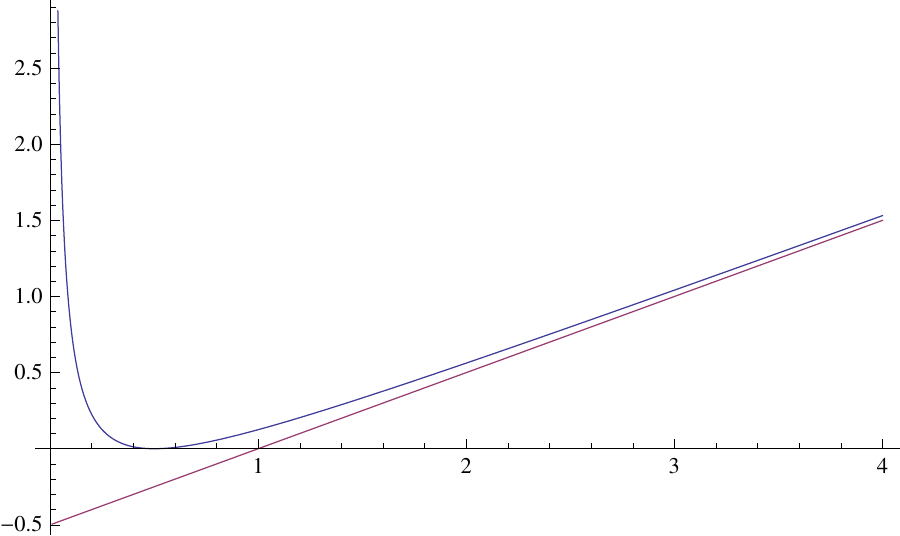}
\end{center}
\caption{Example \ref{ex1}: The rate function  ${\cal I}$ for $\nu=1$.}
\end{figure}
\pagebreak
\par
Let us now comment the change of probability. 
Under $\mathbb P^m$,  the processes $(\xi_t)$ and $(X_t)$ remain in the same family as under $\mathbb P$ : they are Brownian motion with drift and squared Bessel, respectively, but with $\nu$ replaced by $\nu+2m$. 
 Moreover, it is known from \cite{yor2001exponential} Chap. 8 (or equivalently \cite{CPY3}), that under $\mathbb P$,  \[I_\infty \el \left(2Z_\nu\right)^{-1}\] where $Z_\nu$ is a gamma variable of parameter $\nu$. This allows to see directly that the assertion of Lemma \ref{finite} holds.

\subsection{Lévy processes of the form $\xi_t = \pm \D t \pm \hbox{Pois}(\beta, \gamma)_t$}
\label{Poisson}
Here  Pois$(a,b)_t$ is the compound Poisson process of parameter $a$ whose jumps are exponential r.v. of parameter $b$.
They are studied in Section 8.4.3 of \cite{yor2001exponential} (or equivalently \cite{CPY3}), with very informative Tables. In particular, these families are invariant by changes of probability  and the law of the exponential functional $I_\infty$ is known, which allows again  to see directly that the assertion of Lemma \ref{finite} holds.

\subsubsection{$\xi_t =  \D t +  \hbox{Pois}(\beta, \gamma)_t$ with $\D \geq 0$}\label{ex2.1}
(The particular case of the compound Poisson process corresponds to the case $\D= 0$).
We have 
\[\psi(m) = m\left(\D + \frac{\beta}{\gamma-m}\right) , \ m_- = -\infty \ \ , \ m_+ = \gamma \ , \ \psi(m_+) = \infty\,,\]
Here  $m_0 = -\infty$ and $\psi(m_0) = -\infty$ so that
\[\tau_+ = 0\ ,\ \tau_0 = \D^{-1}\ ,\ \tau_e = \frac{\gamma}{\gamma\D + \beta}\ ,\ \Delta = (0,  \D^{-1})\,.\]

The rate function is
\[\mathcal I(x) =  \left(\sqrt{\gamma(1 -  \D x)}- \sqrt{\beta x}\right)^2\ \ , \  \ (0 \leq x \leq \D^{-1})\,.\]

\bigskip\par
\begin{figure}[ht!]
\begin{center}
\includegraphics[height= 160pt]{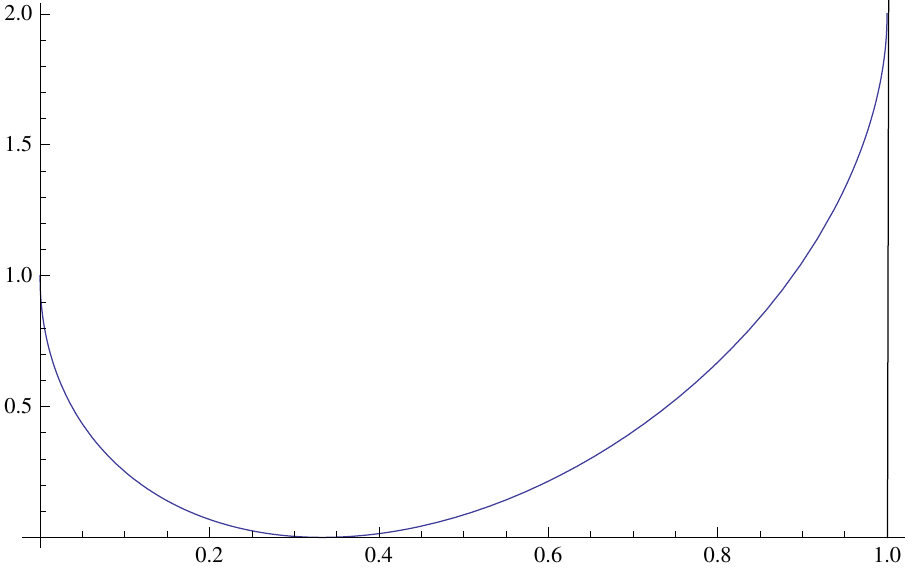}
\end{center}
\caption{Example \ref{ex2.1}: The rate function ${\cal I}$ for 
$\beta=2 , \gamma=1, \D=1$.}
\end{figure}

To check the assumption  (\ref{assum}) we go back to  the process $\xi$. Notice that $\xi_s \geq \D s$ for every $s \geq 0$, so that $\tau^{(\xi)}(t) \leq \D^{-1} \log (1+ \D t)$, whence, for every $\varepsilon > 0$ 
\[\mathbb P \left(\frac{\tau^{(\xi)}(t)}{\log t} > \tau_0 + \varepsilon\right) = 0\]
for $t$ large enough, which is equivalent to (\ref{assum}) thanks to (\ref{fund}).

The boundary $\tau_+$ is in situation 4a and $\tau_0$ in situation 3b, i.e. vertical tangents.

\medskip

In the following examples, with no loss of generality we may and will assume that $\D =1$, since (with obvious notations)
\[\tau^{(\xi ; \D, \beta)}(t) \el \D^{-1}\tau^{(\xi ; 1 , \beta/\D)} (\D \, t)\,.\]
\subsubsection{$\xi_t = -  t +  \hbox{Pois}(\beta, \gamma)_t$ with $0 < \gamma < \beta$}\label{ex2.2}
We have 
\[\psi(m) = m\left(-1 + \frac{\beta}{\gamma-m}\right) , m_- = -\infty \ \ , \ m_+ = \gamma \ , \ \psi(m_+) = \infty\,.\]
Here
\[m_0 = \gamma - \sqrt{\beta\gamma} \ , \ \psi(m_0) = - (\sqrt{\gamma} - \sqrt \beta)^2\,,\]
so that
\[\tau_+ = 0 \ , \ \tau_0 = \infty\ , \ \tau_e = \frac{\gamma}{\beta-\gamma}\ , \ \Delta = (0, \infty)\,.\]
The rate function is
\[\mathcal I(x) = \left(\sqrt{\gamma(1 +  x)}- \sqrt{\beta x}\right)^2\ \ \ (x > 0)\,.\]
The boundary $\tau_+$ is in situation 4a and $\tau_0$ in situation 3a with asymptote $y= (\sqrt \beta - \sqrt \gamma)^2 x  + \gamma - \sqrt{\beta\gamma}$.
\begin{figure}[ht!]
\begin{center}
\includegraphics[height= 160pt]{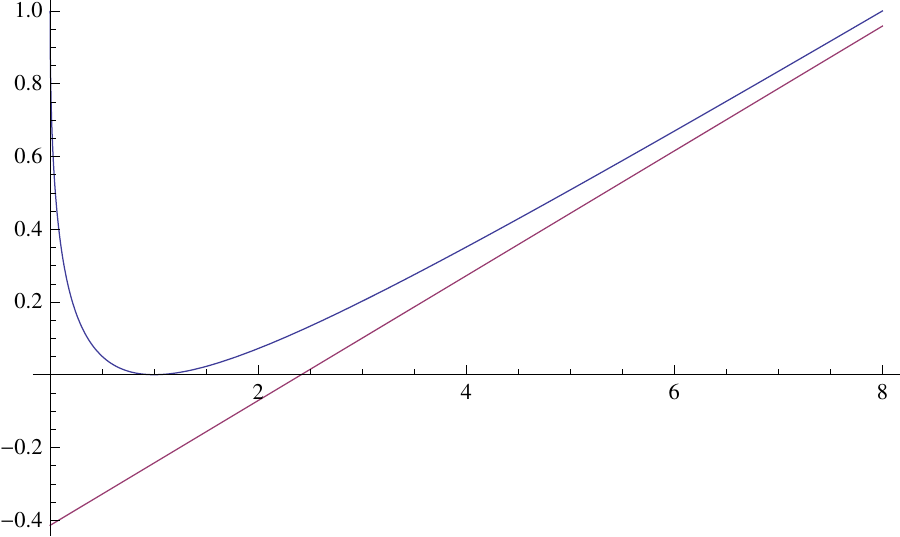}
\end{center}
\caption{Example \ref{ex2.2}: The rate function ${\cal I}$ for 
$\beta=2 , \gamma=1$.}
\end{figure}

\subsubsection
{The saw-tooth process}\label{sawtooth}
\label{saw}

It is a particular case of  the Cramer-Lundberg risk process (\cite{Kyp} Section 1.3.1).
$0 < \beta < \gamma$. 
The Lévy process $\xi$ 
 is defined  by
\begin{equation}
\label{Pois}
\xi
_t = t - \hbox{Pois}
(\beta, \gamma)_t\,.
\end{equation}
The self-similar process $X_t$ is described precisely in (\cite{ab} p. 327). 

 The Lévy exponent of $(\xi_t)$ is   given by:
\begin{equation}
\label{levxi}
\psi
(m) := m\frac{\gamma-\beta+m}{\gamma+m}
  \ \ , \  m_- =-\gamma
 \  , \ m_+ = \infty \ , \ \psi(m_+) = \infty\ , \ \psi'(m_+) = 1\,.\end{equation} 
We have easily
\[m_0 = -
\gamma +\sqrt{\beta\gamma}\ , \ \psi(m_0) = - \left(\sqrt\gamma -\sqrt \beta\right)^2\]
so that 
\[\tau_+ = 1\ , \ \tau_0 = \infty\ , \ \tau_e = \frac{\gamma}{\gamma-\beta}\ , \  \Delta = (1, \infty)\,.\]
The rate function is
\[{\mathcal I}(x) = 
\left(\sqrt{\gamma
(x-1)}- \sqrt{\beta x}\right)^2 \  \ , \ (x \geq 1)\,.\]
The condition (\ref{assum}) is fulfilled since by definition 
  $\xi_t \leq t$ for all $t > 0$, so we have  $\tau^{(\xi)}(t) \geq \log t$ for all $t > 0$.
The minimum of the rate function $\mathcal I$ is reached at 
$x= \frac{
\gamma}{\gamma-\beta}$
 in accordance with the law of large numbers (\cite{ab}, Theorem 4.7 ii)).

The boundary $\tau_+=1$ is in situation 4b and the boundary $\tau_0$ in situation 3a, with
 an asymptote of equation $y= (\sqrt \gamma - \sqrt \beta)^2x +\sqrt{\beta\gamma} - \gamma$.

\bigskip\par
\begin{figure}[ht!]
\begin{center}
\includegraphics[height= 160pt]{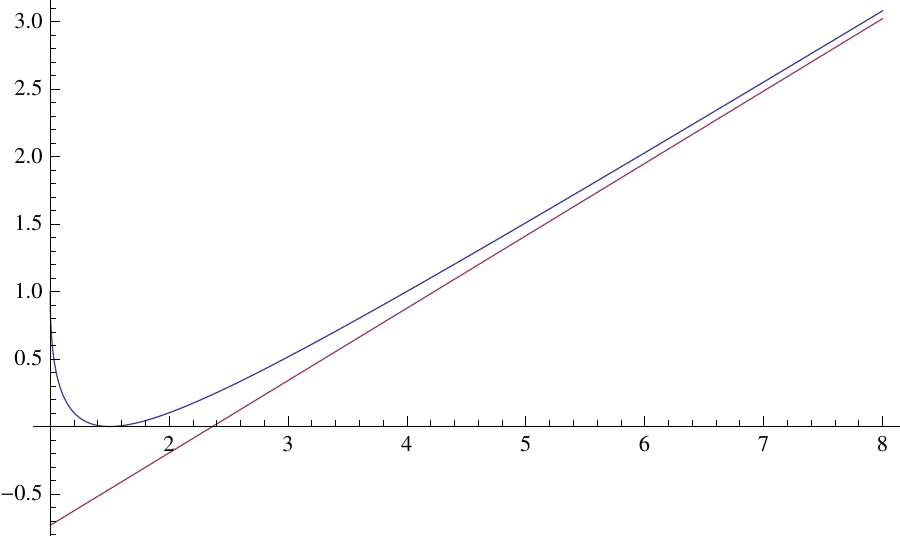}
\end{center}
\caption{Example \ref{sawtooth}: The rate function ${\cal I}$ for 
$\beta=1 , \gamma=3$.}
\end{figure}

\subsection{Three more examples}
\label{threemore}
The first and the second may be seen for instance in \cite{patie1}. The third one is a general class (see \cite{kypext}).
\subsubsection{First}
For $\alpha \in (1,2)$, 
let $X^{\uparrow}$ be the spectrally  negative regular $\alpha$-stable Lévy process conditioned to stay positive. It is a pssMp of index $\alpha$ and its associated Lévy process $\xi^{\uparrow}$ has Laplace exponent
\[\psi(m) = c  \frac{\Gamma(m+\alpha)}{\Gamma(m)}\]
and
\[ m_- = -\alpha \ \ , \ m_+ = \infty\ , \ \psi(m_+) = \infty\ , \ \psi'(m_+) = \infty\]
where $c$ is a positive constant. 
We have
\[\frac{\psi'(m)}{\psi(m)} =  \Psi(m + \alpha) - \Psi(m)\]
where 
 $\Psi$ is the Digamma function.
Then $m_0 = \gamma_\alpha$ where $\gamma_\alpha$ is the (unique) solution in $(-1, 0)$ of the equation
\[\Psi (\gamma+\alpha) = \Psi(\gamma)\,,\]
and then $\psi(\gamma_\alpha) = c\frac{\Gamma(\gamma_\alpha + \alpha)}{\Gamma(\gamma_\alpha)} < 0$.
We have
\[\tau_+ = 0\  , \  \tau_0 = \infty\  , \tau_e = \frac{1}{c\Gamma(\alpha)}\  , \  \Delta = (0, \infty)\,.\]
The boundary $\tau_+$ is in situation 4c and $\tau_0$ in situation 3a with asymptote.

Notice that when $\alpha \uparrow 2$ we end up with the Laplace exponent of a Brownian motion with drift  $1/2$.
\subsubsection{Second}
Let $\kappa \in (0, 1]$ and $\delta > \kappa/(1+\kappa)$. Let $X$ be the continuous state branching process with immigration whose branching mechanism (see \cite{patie1} Lemma 4.8.) is
\[\varphi(u) = - \frac{c}{\kappa} u^{\kappa +1}\]
and immigration mechanism is 
\[\chi(u) = c \delta \frac{\kappa+1}{\kappa} u^{\kappa}\,.\]
It is a pssMp of index $\kappa$ and 
 the associated Lévy process has Laplace exponent
\[\psi (m) = c (\kappa - (\kappa +1) \delta -m)\frac{\Gamma(-m+\kappa)}{\Gamma(-m)}  \ , \ m_- = -\infty \ ,  \ m_+= \kappa \ , \psi(m_+) = \infty\,.\]
First, we see that $\psi'(0) = \lim \psi(m)/m = c ((\kappa+1)\delta- \kappa) \Gamma(\kappa) > 0$ and $\tau_e = 1/ \psi'(0)$. 
Since as $m \rightarrow -\infty$ we have $\psi(m) \sim c(-m)^{\kappa}\rightarrow \infty$, we deduce the existence of $m_0 > -\infty$ such that $\psi'(m_0)= 0$ (and $\psi(m_0) < 0$).
Besides, when $m \uparrow \kappa $ or $h = \kappa-m \downarrow 0$ then 
\[\psi'(m) \sim \frac{(\kappa +1) \delta\Gamma'(h)}{\Gamma(-\kappa)} \uparrow \infty\,.\]

We have then 
\[\tau_+ = 0 \ , \  \tau_0 = \infty \ , \ \Delta = (0, \infty)\,.\]
The boundary $\tau_+=0$ is in situation 4c and $\tau_0$ is in situation 3a with asymptote.

\subsubsection{Hypergeometric-stable process}\label{Cauchy}
The modulus of a Cauchy process in $\mathbb R^d$ for $d > 1$ is a pssMp of index 1 with infinite lifetime. Actually the associated Lévy process is a particular case of hypergeometric-stable process of index $\alpha$ as defined in \cite{cab}, with $\alpha < d$.

The characteristic exponent  is given therein by  Theorem 7, hence the Laplace exponent is :
\[\psi(m) = - 2^\alpha \frac{\Gamma((- m +\alpha)/2)}{\Gamma(-m/2)} \frac{\Gamma((m+d)/2)}{\Gamma((m+d-\alpha)/2)}\]
and
\[  m_- = -d\ , \ m_+=\alpha\ , \ \psi(m_+) = \infty\,.\]

Since $m\Gamma(-m/2) = -2 \Gamma((2-m)/2)$ we have\[\psi'(0)= \lim_{m\rightarrow 0} \psi(m)/m = 2^{\alpha -1}\frac{\Gamma(\alpha/2)\Gamma(d/2)}{\Gamma((d-\alpha)/2)} > 0\,.\]

Moreover, for $m \not= 0$
\[\frac{\psi'(m)}{\psi(m)} = -{\Psi}(((- m +\alpha)/2)) + \Psi((m+d)/2) + \Psi(-m/2) - \Psi(m+d-\alpha)/2)\,,\]
where $\Psi$ is the digamma function.

It is then easy to see that $\psi'$ vanishes at $m = m_0 := (\alpha-d)/2$ and 
\[\psi(m_0) = -2^\alpha\left(\frac{\Gamma((d+ \alpha)/4)}{\Gamma((d- \alpha)/4)}\right)^2\,.\] 
We have then
\[\tau_+ = 0 \ , \ \tau_0 = \infty \ , \ \Delta = (0, \infty)\,.\]

In the particular case $\alpha=1$ and $d=3$, applying the identity
\[\Gamma(x) \Gamma(1-x) = \frac{\pi}{\sin \pi x}\]we have  
\begin{equation}\label{petit}\psi(m) =  (m+1) \tan\frac{\pi m}{2}\end{equation}
with $m_-=-3 \ \ , \ \ m_+ = 1$ and  $\psi(-1) = -2/\pi$ by continuity.
A direct computation gives
\[\psi'(m) 
= \frac{\sin \pi m + \pi m+ \pi}{2 \cos^2\frac{\pi m}{2}}\]
so that 
\[\psi'(0) = \pi/2\,,\quad m_0 = -1\,,\quad \theta_0 = -\psi(m_0) = 2/\pi\mbox{ and }\Delta = (0, \infty)\] We have no close expression for the rate function.
The asymptote as $x \rightarrow \infty$ is the line $y = \frac{2x}{\pi} + 1$. 

\begin{figure}[ht!]
\begin{center}
\includegraphics[height= 160pt]{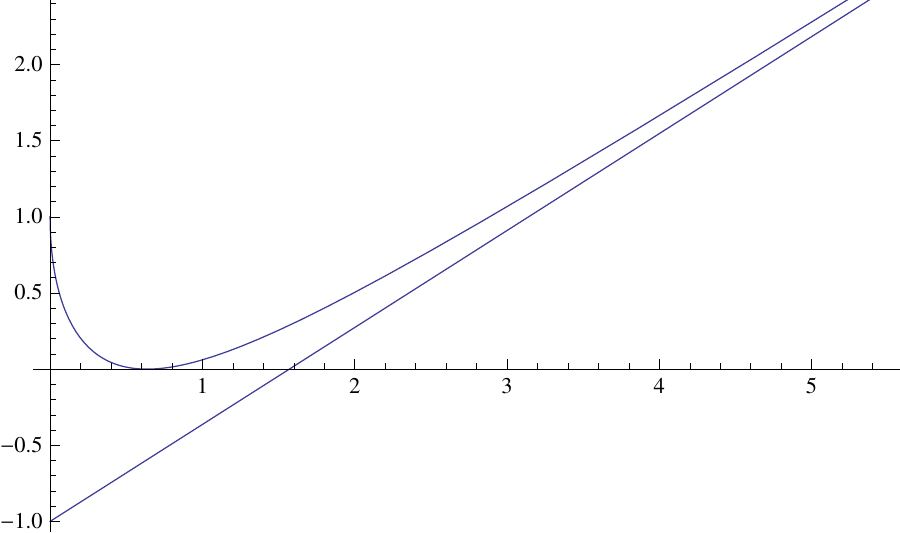}
\end{center}
\caption{Example \ref{Cauchy}: The rate function ${\cal I}(x)$ for $\alpha=1$,
$\D=3$.}
\end{figure}

Notice that in \cite{CPY}, the authors announced a study of the modulus of a multidimensional Cauchy process.  
They found the expression (\ref{petit}) 
 but never published it  (\cite{Pt}). 
\bigskip

\ackname : The authors want  to thank Frédérique Petit for valuable conversations on the Cauchy clock.

\bibliographystyle{plain}

\end{document}